\documentclass[11pt]{article}
\usepackage{amsmath,amssymb,amsthm}
\usepackage{a4wide}
\newtheorem{thm}{Theorem}
\newtheorem{lmm}[thm]{Lemma}
\newtheorem{prp}[thm]{Proposition}

\theoremstyle{definition}

\theoremstyle{remark}
\newtheorem{rem}{Remark}

\newcommand{\R}{\mathbb{R}}
\newcommand{\C}{\mathbb{C}}
\newcommand{\N}{\mathbb{N}}

\newcommand{\la}{\langle}
\newcommand{\ra}{\rangle}
\newcommand{\pd}{\partial}

\newcommand{\bM}{\mathbb{M}}

\DeclareMathOperator{\sgn}{sgn}
\begin{document}
\title{Asymptotic stability of small solitons to 1D NLS with potential}
\author{Tetsu Mizumachi
\thanks{
Faculty of Mathematics Kyushu University Fukuoka 812-8581 Japan,
mizumati@math.kyushu-u.ac.jp}}
\date{}
\maketitle
\begin{abstract}
We consider asymptotic stability of a small solitary wave to
supercritical $1$-dimensional nonlinear Schr\"{o}dinger equations
$$ iu_t+u_{xx}=Vu\pm |u|^{p-1}u \quad\text{for 
$(x,t)\in\mathbb{R}\times\mathbb{R}$,}$$
in the energy class. This problem was studied by
Gustafson-Nakanishi-Tsai \cite{GNT} in the $3$-dimensional case
using the endpoint Strichartz estimate.

To prove asymptotic stability of solitary waves,
we need to show that a dispersive part $v(t,x)$ of a solution belongs to
$L^2_t(0,\infty;X)$ for some space $X$. In the $1$-dimensional
case, this property does not follow from the Strichartz estimate alone.

 In this paper, we prove that a local smoothing estimate
of Kato type holds globally in time and combine the estimate with
the Strichartz estimate to show
$\|(1+x^2)^{-3/4}v\|_{L^\infty_xL^2_t}<\infty$,
which implies the asymptotic stability of a solitary wave.
\end{abstract}

\section{Introduction}
\label{intro}
In this paper, we consider asymptotic stability of solitary wave
solutions to
\begin{equation}
  \label{eq:NLS}
\left\{
  \begin{aligned}
& iu_t+u_{xx}=Vu+f(u) \quad\text{for 
$(x,t)\in\mathbb{R}\times\mathbb{R}$},
\\
& u(x,0)=u_0(x)\quad\text{for $x\in\mathbb{R}$,}
  \end{aligned}\right.
\end{equation}
where $V(x)$ is a real potential, $f(u)=\alpha|u|^{p-1}u$ with $\alpha=\pm1$.

Let 
\begin{align*}
& H(u)=\int_{\R}\left(|u_x|^2+V(x)|u|^2
+\frac{2\alpha}{p+1}|u|^{p+1}\right)dx,  
\\
& N(u)=\int_{\R}|u|^2dx.
\end{align*}
Then a solution to \eqref{eq:NLS} satisfies
\begin{equation}
  \label{eq:cons}
H(u(t))=H(u_0),\quad N(u(t))=N(u_0)  
\end{equation}
during the time interval of existence. Stability of solitary waves to NLS was
first studied by Cazenave and Lions \cite{CaLi},
Grillakis-Shatah-Strauss \cite{GSS} and Weinstein \cite{W2} 
(see also Rose-Weinstein \cite{RW}, Oh \cite{Oh} and Shatah-Strauss
\cite{ShSt}). In the case of integrable equations such as cubic NLS
and KdV, the inverse scattering theory tells us that if the initial data
decays rapidly as $x\to\pm\infty$, a solution decomposes into
a sum of solitary waves and a radiation part as $t\to\infty$
(see \cite{Sc}). Soffer and Weinstein \cite{SW1,SW2}
considered NLS with potential
\begin{equation}
  \label{eq:NLS2}
iu_t+\Delta u=Vu\pm |u|^{p-1}u \quad\text{for $x\in\R^n$ and $t>0$,}  
\end{equation}
where $n\ge 2$ and $1<p<(n+2)/(n-2)$. They proved that if $-\Delta +V$
has exactly one negative eigenvalue and initial data
is well localized and close to a nonlinear bound state, a solution tends
to a sum of a nonlinear bound state nearby and a radiation part which
disperses to $0$ as $t\to \infty$.
This result was extended by Yau and Tsai \cite{YT1,YT2,YT3} and 
Soffer-Weinstein \cite{SW3} to the case where $-\Delta+V$ have two bound
states.
In the $1$-dimensional case, Buslaev and Perelman \cite{BP1, BP2}
and Buslaev and Sulem \cite{BuS} studied the asymptotic stability of
\eqref{eq:NLS} with $V\equiv0$.
Using the Jost functions, they built a local energy decay estimate of
solutions to the linearized equation and prove asymptotic stability 
of solitary waves for super critical nonlinearities (see also \cite{ZGS}). 
Their results are extended to the higher dimensional case by
Cuccagna \cite{Cu} (see also Perelman \cite{Per} and 
Rodnianski-Schlag-Soffer \cite{RSS} which study asymptotic stability
of multi-solitons).
\par
However, all these results assume that initial data is well localized
 so that a solution decays like $t^{-3/2}$.
Martel and Merle \cite{MM}, \cite{MM2} proved the asymptotic stability
of solitary waves to generalized KdV equations using the monotonicity of $L^2$-mass, which is a variant of the local smoothing effect proved by
Kato \cite{K}. They elegantly use the fact that a dispersive remainder
part of a solution $v(t,x)$ satisfies
\begin{equation}
  \label{eq:disp}
\int_0^\infty \|v(t,\cdot)\|_{H^1_{loc}}^2dt<\infty  
\end{equation}
to prove the asymptotic stability of solitary waves in $H^1$
(see also El-Dika \cite{ElD} and Mizumachi \cite{Mi2} for BBM equation
and Pego-Weinstein \cite{PW2} and Mizumachi \cite{Mi1} for
KdV with localized initial data).
Recently, Gustafson-Nakanishi-Tsai \cite{GNT} has proved asymptotic
stability of a small solitary wave of \eqref{eq:NLS2} in the energy
class with $n=3$. Their idea is to use the endpoint Strichartz
estimate instead of \eqref{eq:disp}, which tells us that
$\|v\|_{L^2_tW^{1,6}_x}$ remains small globally in time
for super critical nonlinearity.
However, the Strichartz estimate is not sufficient in the lower 
dimensional case to obtain some estimate like \eqref{eq:disp} because a
dispersive wave decays more slowly than the $3$-dimensional case.
To overcome this difficulty, we prove
\begin{align}
\label{eq:l1}
& \|\la x\ra^{-3/2}e^{it(-\pd_x^2+V)}Qf\|_{L_x^\infty L^2_t}\le
 C\|f\|_{L^2},\\
\label{eq:l2}
& \|\pd_xe^{{it(-\pd_x^2+V)}}Qf\|_{L_x^\infty L^2_t}
\le C\|f\|_{H^{1/2}},
\end{align}
where $Q$ is a spectral projection associated to the continuous spectrum of $-\pd_x^2+V$. The local smoothing estimate of $1/2$ gain derivative 
has been studied by many authors (see e.g. Constantin and Saut \cite{CS},
Kato and Yajima \cite{KY} and Kenig-Ponce-Vega \cite{KPV1,KPV2})
to show the local well-posedness of semilinear equations with derivative terms. Most of them are without potential (\cite{CS,Sj}) or local in
time (see \cite{RZ1}).

Ben-Artzi and Klainerman \cite{Be-Kl} proved a time global local
smoothing estimate for the $n$-dimensional case with $n\ge3$.
 See also Barcel\'{o}-Ruiz-Vega \cite{BRVe} who use
a Morawetz type inequality to obtain the result.
Recently, Burq and Planchon \cite{Bu-Pl} has proved local smoothing
estimates including an estimate similar to \eqref{eq:l2}
for $Lu=-\pd_x(a(x)\pd_xu)$ (they use $\dot{B}^{1/2}_{2,\infty}$ instead
of $L^\infty$).
In the present paper, we show \eqref{eq:l1} and \eqref{eq:l2} assuming
the non-resonance condition for $L=-\Delta+V$.
Another difference between \cite{Bu-Pl} is that $L$ may have negative
eigenvalues.

Our proof given in this paper for the $1$-dimensional case is different
from \cite{Be-Kl,Bu-Pl,BRVe}. We use the Born series
(see Artbazar-Yajima \cite{AY} and Goldberg-Schlag \cite{GSc})
for the high frequency part and a theory
of Jost functions for the low frequency part.

\par Finally, we introduce several notations.
For complex valued functions $f(x)$ and $g(x)$, we denote
$\la f,g\ra=\int_\R f(x)\overline{g(x)}dx$.
Let
\begin{gather*}
\|u\|_{L_t^qL_x^p}=\biggl(\int_\R\bigl(\int_\R |f(t,x)|^pdx\bigr)^{q/p}dt
\biggr)^{1/q},
\\
\|u\|_{L_x^sL_t^r}=\biggl(\int_\R\bigl(\int_\R |f(t,x)|^rdt\bigr)^{s/r}dx
\biggr)^{1/s},
\end{gather*}
and let $H^{1,k}(\R)$ be the Hilbert space equipped with the norm
$$\|u\|_{H^{1,k}}=
\left(\sum_{i=0,1}\int_\R (1+x^2)^k|\pd_x^iu(x)|^2dx\right)^{1/2}.$$
For any Banach spaces $X$, $Y$, we denote by $B(X,Y)$ the space of
bounded linear operators from $X$ to $Y$. We abbreviate
$B(X,X)$ as $B(X)$.

We define the Fourier and transform of $f(x)$ as
$$
\mathcal{F}_xf(\xi)=\hat{f}(\xi)=
\frac{1}{\sqrt{2\pi}}\int_\R f(x)e^{-ix\xi}dx,$$
and the inverse Fourier transform of $g(\xi)$ as
$$
\mathcal{F}^{-1}_\xi g(x)=\hat{g}(-x)
=\frac{1}{\sqrt{2\pi}}\int_\R g(\xi)e^{ix\xi}d\xi.
$$
We define $\mathcal{S}_\otimes(\R^2)$  as a set of functions written as 
$f(t,x)=\sum_{i=1}^N f_i(t)g_i(x)$ with $f_i$, $g_i\in \mathcal{S}(\R)$
($1\le i\le N$).

For an interval $I\subset \R$, let $\chi_I(x)$ be a
characteristic functions satisfying $\chi_I(x)=1$ for $x\in I$ and
$\chi_I(x)=0$ for $x\not\in I$. We denote $\sqrt{1+|x|^2}$ by $\la x\ra$.

\section{The Main result and Preliminaries}
In the present paper, we assume that the linear potential $V(x)$ is a
continuous function on $\R$ and satisfies the following.
\begin{itemize}
\item[(V1)] $(1+x^2)V(x)\in L^1(\R)$.
\item [(V2)] $L=-\pd_x^2+V$ has exactly one negative eigenvalue $E_*$,
and $0$ is neither a resonance nor an eigenvalue of $L$.
\end{itemize}
Let $E\in\R$ and $e^{-iEt}\phi_E(x)$ be a solitary wave solution of
\eqref{eq:NLS}.
Then $\phi_E(x)$ is a solution to
\begin{equation}
\label{eq:B}
\left\{\begin{aligned}
&\phi_E''+E\phi_E=V\phi_E+\alpha|\phi_E|^{p-1}\phi_E 
\quad\text{for $x\in \R$},\\
& \lim_{x\to\pm\infty}\phi_E(x)=0.
\end{aligned}\right.
\end{equation}

Using the bifurcation theory, we have the following.
\begin{prp}
  \label{prop:1} Assume (V1) and (V2). Then there exists a $\delta>0$ such
that if $E\in(E_*,E_*+\delta)$ and $\alpha=1$ or $E\in (E_*-\delta,E_*)$ and
$\alpha=-1$, Eq. \eqref{eq:B} has a positive solution $\phi_E$ satisfying
the following:
\begin{enumerate}
\item $\phi_E\in H^{1,k}$ for every $k\in\N$,
\item the function $E\mapsto \phi_E$ is $C^2$ in
$H^{1,k}$ for every $k\in\N$, and as $E\to E_*$,
$$
\phi_E=|E-E_*|^{1/(p-1)}\left(\|\phi_*\|_{L^{p+1}}^{-(p+1)/(p-1)}
\phi_*+O(E-E_*)\right)\quad\text{in $H^{1,k}$},$$
where $\phi_*$ is an normalized eigenfunction of $H$ {\rm(}satisfying
$\|\phi_*\|_{L^2}=1${\rm)} belonging to $E_*$.
\end{enumerate}
\end{prp}
Proposition \ref{prop:1} follows from a rather standard argument.
See for example \cite{Ni} and \cite[pp.123--124]{SW1}.  
\par

Now, we introduce our main result.
\begin{thm}
\label{thm1} Assume (V1) and (V2).
Let $p\ge5$ and $\varepsilon_0$ be a sufficiently small positive number.
Suppose $\|u_0\|_{H^1}<\varepsilon_0$. Then there exist an $E_+<0$,
a $C^1$ real-valued function $\theta(t)$ and
$v_+\in H^1(\R)$ such that
\begin{gather*}
\lim_{t\to\infty}\dot{\theta}(t)=E_+,\\
|E_+-E_*|+\|v_+\|_{H^1}=O(\|u_0\|_{H^1}),\\
\lim_{t\to\infty}\|u(t)-e^{i\theta(t)}\phi_{E_+}
{}-We^{it\pd_x^2}v_+\|_{H^1(\R)}=0,
\end{gather*}
where $W=\lim_{t\to\infty}e^{-itL}e^{-it\pd_x^2}$.
\end{thm}
\begin{rem}
\label{rem:1}
Let
$\phi_{1,E}=\|\phi_E\|_{L^2}^{-1}\phi_E$ and
$\phi_{2,E}=\|\pd_E\phi_E\|_{L^2}^{-1}\pd_E\phi_E$.
By Proposition \ref{prop:1},
\begin{equation*}
\|\phi_{1,E}-\phi_*\|_{H^{1,k}(\R)}
+\|\phi_{2,E}-\phi_*\|_{H^{1,k}(\R)}\lesssim |E-E_*|.
\end{equation*}
\end{rem}
\begin{rem}  
Let us decompose a solution to \eqref{eq:NLS} into a solitary wave part and
a radiation part:
\begin{equation}
  \label{eq:1.3}
  u(t,x)=e^{-i\theta(t)}(\phi_{E(t)}(x)+v(t)).
\end{equation}
If we take initial data in the energy class, the dispersive part of
the solutions decays more slowly than they does for well localized
initial data. So, being different from Soffer-Weinstein \cite{SW1,SW2}
or Buslaev-Perelman \cite{BP1}, we cannot expect that
$\int_t^\infty\dot{E}(s)ds$ is integrable.
Thus in general, we need dispersive estimates for a time-dependent
linearized equations to prove asymptotic stability of solitary waves
in $H^1(\R)$. To avoid this difficulty, we assume the smallness of
solitary waves so that a generalized kernel of the linearized operator
is well approximated by a $1$-dimensional subspace
$\{\beta \phi_*\,|\,\beta\in\C\}.$
\end{rem}
Substituting \eqref{eq:1.3} into \eqref{eq:NLS}, we obtain
\begin{equation}
  \label{eq:v}
  iv_t=Lv+g_1+g_2+g_3+g_4,
\end{equation}
where 
\begin{align*}
& g_1(t)=-\dot{\theta}(t)v(t), \quad
 g_2(t)=(E(t)-\dot{\theta}(t))\phi_{E(t)}-i\dot{E}(t)\pd_E\phi_{E(t)},
\\ 
& g_3(t)=f(\phi_{E(t)}+v(t))-f(\phi_{E(t)})-\pd_\varepsilon
 f(\phi_{E(t)}+\varepsilon v(t))|_{\varepsilon=0},
\\
& g_4(t)=\pd_\varepsilon f(\phi_{E(t)}+\varepsilon v(t))|_{\varepsilon=0}
=\alpha\phi_{E(t)}^{p-1}\left(\frac{p+1}{2}v(t)
+\frac{p-1}{2}\overline{v(t)}
\right).
\end{align*}

\label{sec:2}
To fix the decomposition \eqref{eq:1.3}, we assume
\begin{equation}
  \label{eq:2.1}
\left\la \Re v(t), \phi_{E(t)}\right\ra
=\left\la \Im v(t),\pd_E\phi_{E(t)}\right\ra=0.
\end{equation}
By Proposition \ref{prop:1}, we have
\begin{equation}
  \label{eq:init-err}
|E(0)-E_*|^{1/(p-1)}+\|v(0)\|_{H^1} \lesssim \|u_0\|_{H^1}.  
\end{equation}
Since $u\in C(\R;H^1(\R))$, it follows from
the implicit function theorem that there exist a $T>0$ and $E$,
$\theta\in C^1([-T,T])$ such that \eqref{eq:2.1} holds for $t\in[-T,T]$.
See, for example, \cite{GNT} for the proof.

Differentiating \eqref{eq:2.1} with respect to $t$ and substituting
\eqref{eq:v} into the resulting equation, we obtain
\begin{equation}
  \label{eq:2.2}
\mathcal{A}(t)
\begin{pmatrix}  \dot{E}(t)\\ \dot{\theta}(t)-E(t)\end{pmatrix}
=\begin{pmatrix}
\la \Im g_3(t),\phi_{E(t)}\ra \\  \la \Re g_3(t),\pd_E\phi_{E(t)}\ra
\end{pmatrix},\end{equation}
where
\begin{align*}
& \mathcal{A}(t)=\\
& \begin{pmatrix}
 \la\pd_E\phi_{E(t)},\phi_{E(t)}\ra-\la \Re v(t),\pd_E\phi_{E(t)}\ra
& \la \Im v(t),\phi_{E(t)}\ra\\
 \la \Im v(t),\pd_E^2\phi_{E(t)}\ra
& \la\pd_E\phi_{E(t)},\phi_{E(t)}\ra+\la \Re v(t),\pd_E\phi_{E(t)}\ra
\end{pmatrix}.  
\end{align*}
\par

To prove our main result, we will use the Strichartz estimate and
a time global estimate of Kato type.
The Strichartz estimate along with $L^\infty-L^1$-estimate for
$1$-dimensional Schr\"{o}dinger equations with linear potential was
obtained by Goldberg and Schlag \cite{GSc}.

Let $Pu=\left\la u, \phi_*\right\ra\phi_*$, $Qu=(I-P)u.$
Then we have the following.
\begin{lmm}[Strichartz estimate (\cite{GSc,KrSc})]
  \label{lem:2.1} Assume (V1) and (V2).
  \begin{itemize}
  \item [(a)]
There exists a positive number $C$ such that for any $f\in L^2(\R)$,
$$\|e^{-itL}Qf\|_{L_t^4L_x^\infty\cap L_t^\infty L_x^2}\le C\|f\|_{L^2}.$$
\item[(b)]
There exists a positive number $C$ such  that for any
$g(t,x)\in \mathcal{S}(\R^2)$,
\begin{gather*}
\|\int_{0}^te^{-i(t-s)L}Qg(s,\cdot)ds\|_{L_t^4L_x^\infty\cap L_t^\infty L_x^2}
\le C\|g\|_{L_t^{4/3}L_x^1+L_t^1L_x^2}.
\end{gather*}
  \end{itemize}
\end{lmm}
To estimate the the quadratic term of $v$ in $g_3$,
we need the following lemma.
\begin{lmm}
  \label{lem:2.2}
Assume (V1) and (V2).
  \begin{enumerate}
  \item [(a)]
There exists a positive constant $C$ such that for any
$f\in \mathcal{S}(\R)$,
\begin{align}
\label{eq:2.2.1}
& \|\la x\ra^{-3/2}e^{-itL}Qf\|_{L_x^\infty L^2_t}\le
 C\|f\|_{L^2},\\
\label{eq:2.2.2}
& \|\pd_xe^{-itL}Qf\|_{L_x^\infty L^2_t}
\le C\|f\|_{H^{1/2}}.
\end{align}
  \item [(b)]
There exists a positive constant $C$ such that for any
$g(t,x)\in \mathcal{S}(\R^2)$,
\begin{equation}
\label{eq:2.2.3}
\left\|\int_\R e^{isL}Qg(s,\cdot)ds\right\|_{L^2_x}
\le C\|\la x\ra^{3/2}g\|_{L_x^1L_t^2},
\end{equation}
\end{enumerate}
\end{lmm}
\begin{lmm}
  \label{lem:2.3}
There exists a positive constant $C$ such that for any
$g(t,x)\linebreak\in \mathcal{S}(\R^2)$ and $t\in\R$,
\begin{equation}
\label{eq:2.3.1}
\sum_{j=0,1}\left\|\la x\ra^{-1}\pd_x^j \int_0^t
e^{-i(t-s)L}Qg(s,\cdot)ds\right\|_{L^\infty_xL_t^2}
\le C\|\la x\ra g\|_{L_x^1L_t^2}.
\end{equation}
Furthermore, if $\sup_{x\in\R}e^{\alpha|x|}|V(x)|<\infty$ holds
for an $\alpha>0$, there exists a positive number $C$ such that
\begin{equation}
\label{eq:2.3.2}
\left\|\int_0^t \pd_x e^{-i(t-s)L}Qg(s,\cdot)ds
\right\|_{L^\infty_xL_t^2}\le C\|g\|_{L_x^1L_t^2}.
\end{equation}
\end{lmm}
Lemma \ref{lem:2.1} is not applicable to a linear term $g_4$ in \eqref{eq:v}
because we do not have $g_4\in L^{4/3}_tL_x^1+L^1_tL_x^2$.
To deal with $g_4$, we use a lemma by Christ and Kiselev \cite{ChK} to
combine Lemmas \ref{lem:2.1} and \ref{lem:2.2}.
\begin{lmm}
  \label{lem:2.4} Assume (V1) and (V2).
Then there exists a positive constant $C$ such that for any
$g(t,x)\in \mathcal{S}(\R^2)$  and $t\in\R$,
$$
\left\|\int_0^t e^{-i(t-s)L}Qg(s,\cdot)ds
\right\|_{L_t^4L_x^\infty\cap L_t^\infty L_x^2}
\le C\|g\|_{L_t^2L_x^2(\R;\la x\ra^5dx)}.$$
\end{lmm}
The proof of Lemmas \ref{lem:2.2}--\ref{lem:2.4} will be given 
in Section 4.
\section{Proof of Theorem \ref{thm1}}
\label{sec:3}
In this section, we will prove Theorem \ref{thm1}.
To eliminate $g_1$ in \eqref{eq:v}, we put
\begin{equation}
  \label{eq:trans}
w(t)=e^{-i\theta(t)}v(t).
\end{equation}
Then  \eqref{eq:v} is translated into the integral equation
\begin{equation}
  \label{eq:v'}
w(t)=e^{-itL}w(0)-i\sum_{2\le j\le 4}
\int_0^t e^{-i(t-s)L} e^{-i\theta(s)}g_j(s)ds.
\end{equation}
All nonlinear terms in \eqref{eq:v'} can be estimated in terms of
the following. 
\begin{align*}
& \bM_1(T)=\sup_{0\le t \le T}|E(t)-E_*|,\quad
\bM_2(T)=\|\la x\ra^{-3/2}Qw\|_{L^\infty_xL^2(0,T)},
\\ &
\bM_3(T)=\|Pw\|_{L^\infty_xL^2(0,T)}+\|\pd_x Pw\|_{L^\infty_xL^2(0,T)},
\\ & 
\bM_4(T)=\|Qw\|_{L^q(0,T;W^{1,2p}_x)\cap L^\infty(0,T;H^1_x)}
+\|Qw\|_{L^4(0,T;L^\infty_x)},
\\ &
\bM_5(T)=\|Pw\|_{L^4(0,T;W^{1,\infty}_x)\cap L^\infty(0,T;H^1_x)},
\quad \bM_6(T)=\|\pd_x Qw\|_{L^\infty_xL^2(0,T)}.
\end{align*}
where $4/q=1-1/p$.
\begin{proof}[Proof of Theorem \ref{thm1}]
Proposition \ref{prop:1}, Remark \ref{rem:1}  and \eqref{eq:2.1} imply that
$$\la \pd_E\phi_E,\phi_E\ra=O(|E-E_*|^{2/(p-1)-1}),\quad
|\la v, \pd_E^i\phi_E\ra|\lesssim |E-E_*|^{p/(p-1)-i}\|v\|_{L^2}.$$
Thus by \eqref{eq:2.2}, we have
\begin{align}
\label{eq:3.1}
&|\dot{\theta}(t)-E(t)|
 \lesssim \|\phi_{2,E(t)}v^2\|_{L^1}+\|\phi_{2,E(t)}f(v)\|_{L^1},
\\
\label{eq:3.2}
& |\dot{E}(t)|\lesssim \|\phi_{1,E(t)}v^2\|_{L^1}
+\|\phi_{1,E(t)}f(v)\|_{L^1}.
\end{align}

Suppose that the decomposition \eqref{eq:1.3} with
\eqref{eq:2.1} persists for $0\le t\le T$ and that
$\bM_i(T)$ $(1\le i \le 8)$ are bounded.
Eqs. \eqref{eq:3.1}--\eqref{eq:3.2} imply that
\begin{equation}
  \label{eq:3.3a}
  \begin{split}
& \|\dot{\theta}-E\|_{L^1(0,T)}+\|\dot{E}\|_{L^1(0,T)}
\\ \le &
C(\bM)(\|\phi_{1,E(t)}v^2\|_{L^1(0,T;L^1_x)}
+\|\phi_{2,E(t)}v^2\|_{L^1(0,T;L^1_x)})
\\ \le &
C(\bM)\left(\sum_{i=1,2}\left\|\la x\ra^{3}\phi_{i,E(t)}
\right\|_{L^1_xL^\infty(0,T)}\right)
\|\la x\ra^{-3/2}v\|_{L_x^\infty L^2(0,T)}^2
\\ \le &
C(\bM)\left(\sum_{i=1,2}\left\|\la x\ra^{5}\phi_{i,E(t)}
\right\|_{L^\infty(0,T;L^\infty_x)}\right)
\|\la x\ra^{-3/2}v\|_{L_x^\infty L^2(0,T)}^2
\\ \le & C(\bM)(\bM_2(T)+\bM_3(T))^2,
  \end{split}
\end{equation}
and
\begin{equation}
  \label{eq:3.3b}
  \begin{split}
\|\dot{\theta}-E\|_{L^\infty(0,T)}+\|\dot{E}\|_{L^\infty(0,T)}
\lesssim & \sup_{0\le t\le T}(\|v\|_{H^1}^2+\|v\|_{H^1}^p)
\\ \le & C(\bM)(\bM_4(T)+\bM_5(T))^2.
  \end{split}
\end{equation}
Hereafter we denote by $C(\bM)$ various functions of
$\bM_1$, \dots, $\bM_4$ that are bounded in a finite neighborhood of $0$.
By \eqref{eq:3.3a} and \eqref{eq:init-err},
\begin{equation}
  \label{eq:3.4}
\bM_1(T) \lesssim \|u_0\|_{H^1}+ C(\bM)(\bM_2+\bM_3)^2.
\end{equation}

By Remark \ref{rem:1} and \eqref{eq:2.1}, we have
\begin{align*}
  |\la w(t), \phi_*\ra|\le & \|\la x\ra^{-5/2}v\|_{L^2_x(\R)}
\sum_{i=1,2}\|\la x\ra^{5/2}(\phi_{i,E}-\phi_*)\|_{L^2(\R)}
\\ \lesssim & |E(t)-E_*|\|\la x\ra^{-5/2}w\|_{L^2_x(\R)},
\end{align*}
and that
\begin{equation}
  \label{eq:3.6}
  \begin{split}
\bM_3(T) \le & \sup_{t\in[0,T]}|E(t)-E_*|\|\la x\ra^{-5/2}w\|_{L^2(0,T;L_x^2)}
\\ \le & C(\bM)\bM_1(T)(\bM_2(T)+\bM_3(T)).
  \end{split}
  \end{equation}
Similarly, we have
\begin{equation}
    \label{eq:3.7}
\bM_5(T) \le C(\bM)\bM_1(T)(\bM_4(T)+\bM_5(T)).      
\end{equation}

Next, we will estimate $\bM_2(T)$.
By \eqref{eq:v'},
$$\bM_2(T)\le I_1+I_2+I_3+I_4,$$
where
$$
 I_1=\|\la x\ra^{-3/2} e^{-it L}Qw(0)\|_{L_x^\infty L^2(0,T)}
\lesssim  \|w(0)\|_{L^2},
$$
and
$$
I_i=\left\|\la x\ra^{-3/2} \int_0^t e^{-i(t-s)L}Qg_i(s)ds
\right\|_{L^\infty_xL^2(0,T)} \quad\text{for $2\le i\le 4$.}
$$
By Lemma \ref{lem:2.3}, Remark \ref{rem:1}, \eqref{eq:3.3a} and
\eqref{eq:3.3b},
\begin{align*}
  I_2  \lesssim & \|\la x\ra Qg_2\|_{L_x^1L^2(0,T)}
\\  \le &
\left\|\la x\ra Q\phi_{E(t)}\right\|_{L^1_xL^\infty(0,T)}
\|\dot{\theta}-E\|_{L^2(0,T)}
\\ &
+\left\|\la x\ra Q\pd_E\phi_{E(t)}\right\|_{L^1_xL^\infty(0,T)}
\|\dot{E}\|_{L^2(0,T)}
\\  \le &
C(\bM)\bM_1(T)^{1/(p-1)}(\bM_2(T)+\bM_3(T)+\bM_4(T)+\bM_5(T))^2.
\end{align*}
Lemmas \ref{lem:2.2} and \ref{lem:2.3} yield
\begin{align*}
  I_3 \lesssim & \|\la x\ra \phi_{E(t)}^{p-2}v^2\|_{L_x^1L^2(0,T)}
+\int_0^T\||v|^p\|_{L^2_x}ds
\\ \lesssim & \left\|\la x\ra^{5/2}\sup_{t\in[0,T]}\phi_{E(t)}^{p-2}
\right\|_{L^1_x}
\|\la x\ra^{-3/2}v\|_{L_x^\infty L^2(0,T)}
\|v\|_{L^\infty(0,T;L_x^\infty)} 
\\ & 
+ \|v(t)\|_{L^\infty(0,T;L^\infty_x)}^{p-q}
\|v(t)\|_{L^q(0,T;L_x^{2p})}^q
\\ \le & C(\bM)\sum_{2\le i \le 5} \bM_i(T)^2.
\end{align*}
where $4/q=1-1/p$. 
By Proposition \ref{prop:1} and Lemma \ref{lem:2.3},
\begin{align*}
  I_4\lesssim  \|\la x\ra g_4\|_{L_x^1L^2(0,T)}
 \lesssim & \left\|\la x\ra^{5/2} \sup_{t\in[0,T]} \phi_{E(t)}^{p-1}
\right\|_{L^1_x}\|\la x\ra^{-3/2}v\|_{L_x^\infty L^2(0,T)}
\\ \le & C(\bM)\bM_1(T)(\bM_2(T)+\bM_3(T)).
\end{align*}
Combining the above, we see that
\begin{equation}
  \label{eq:3.8}
  \bM_2(T)\le \|u_0\|_{H^1}+C(\bM)\sum_{1\le i\le 5}\bM_i(T)^2.
\end{equation}
Likewise, we have
\begin{equation}
  \label{eq:3.9}
  \bM_6(T)\le \|u_0\|_{H^1}+C(\bM)\sum_{1\le i\le 5}\bM_i(T)^2.
\end{equation}
Finally, we will estimate $\bM_4(T)$.
In view of \eqref{eq:v'},
$$\bM_4(T)\le J_1+J_2+J_3+J_4,$$
where
\begin{align*}
J_1=&\left\| e^{-itL}Qw(0)
\right\|_{L^\infty(0,T;H^1_x)\cap L^q(0,T;W_x^{1,2p})
\cap L^4(0,T;L_x^\infty)}
\\  J_i=&\left\|\int_0^t e^{-i(t-s)L}Qg_i(s)ds
\right\|_{L^\infty(0,T;H^1_x)\cap L^q(0,T;W_x^{1,2p})\cap L^4(0,T;L_x^\infty)}
\end{align*}
for $2\le i\le 4$, where $4/q=1-1/p$.
By Lemma \ref{lem:2.1} and \eqref{eq:3.3a},
$$  J_1 \lesssim \|w(0)\|_{H^1},$$
and
\begin{align*}
J_2 \lesssim & \|Qg_2(s)\|_{L^1(0,T;H^1_x)}
\\ \lesssim &
\|\dot{\theta}-E\|_{L^1(0,T)}\sup_{t\in[0,T]}\|Q\phi_{E(t)}\|_{H^1_x}
+\|\dot{E}\|_{L^1(0,T)}\sup_{t\in[0,T]}\|Q\pd_E\phi_{E(t)}\|_{H^1_x}
\\ \le & C(\bM)(\bM_2(T)^2+\bM_3(T)^2).
\end{align*}
Note that $\|Q\pd_E\phi_E\|_{H^1}\lesssim |E-E_*|^{1/(p-1)}$
follows from Remark \ref{rem:1}.
Using Minkowski's inequality and Lemma \ref{lem:2.1}, we have
\begin{align*}
 J_3 \lesssim  & \|\phi_{E(t)}^{p-2}v^2\|_{L^{6/5}(0,T;W_x^{1,6/5})}
+ \|f(v)\|_{L^1(0,T;H^1_x)}
\\ \lesssim &
\|\la x\ra^2 \mathstrut{\sup_{t\in[0,T]}}\phi_{E(t)}^{p-2}
\|_{W_x^{1,6/5}} \|v\|_{L^4(0,T;L_x^\infty)}^{2/3}
\left(\sum_{i=0,1}\|\la x\ra^{-3/2}\pd_x^iv\|_{L_x^\infty L^2(0,T)}
\right)^{4/3}
\\ & + \|v(t)\|_{L^\infty(0,T;H^1_x)}^{p-q}
\|v(t)\|_{L^q(0,T;W_x^{1,2p})}^q
\\ \le &
C(\bM)\sum_{2\le i\le 6} \bM_i(T)^2,
\end{align*}
where $4/q=1-1/p$.
By Proposition \ref{prop:1} and Lemma \ref{lem:2.4},
\begin{align*}
J_4 \lesssim & \|\la x\ra^{5/2} g_4\|_{L^2(0,T;L_x^2)}
\\ \lesssim & \|\sup_{t\in[0,T]}\la x\ra^4\phi_{E(t)}^{p-1}\|_{L_x^2}
\|\la x\ra^{-3/2}v\|_{L_x^\infty L^2(0,T)}
\\ \lesssim & \bM_1(T)(\bM_2(T)+\bM_3(T)).  
\end{align*}
Combining the above, we have
\begin{equation}
  \label{eq:3.10}
  \bM_4(T)\le \|u_0\|_{H^1}+C(\bM)\sum_{1\le i\le 6}\bM_i(T)^2.
\end{equation}
It follows from \eqref{eq:3.4}--\eqref{eq:3.10} that
if $\varepsilon_0$ is sufficiently small,
\begin{equation}
  \label{eq:3.12}
\sum_{1\le i\le 6}\bM_i(T)\lesssim \|u_0\|_{H^1}.  
\end{equation}
Thus by continuation argument, we may let $T\to\infty$.

By \eqref{eq:3.3a}, there exists an $E_+<0$ satisfying
$\lim_{t\to\infty}E(t)=E_+$ and $|E_+-E_*|\lesssim \|u_0\|_{H^1}.$
In view of \eqref{eq:3.12}, we have
\begin{align*}
w_1:=&-i\sum_{2\le j\le 4}\int_0^\infty e^{isL}Qe^{-i\theta(s)}g_j(s)ds
\in H^1(\R),\\
\|w_1\|_{H^1}\lesssim & \|g_2(s)\|_{L^1_tH^1_x}+
\|g_3\|_{L_t^{6/5}W_x^{1,6/5}+ L_t^1H_x^1}
+\|\la x\ra^{5/2} g_4\|_{L_t^2H_x^1}
\\ \lesssim & \|u_0\|_{H^1},
\end{align*}
and $$\lim_{t\to\infty}\|Qw(t)-e^{-itL}(Qw(0)+w_1)\|_{H^1}=0.$$
By \cite{GSc}, we have
$\|e^{-itL}Qf\|_{L^\infty}\lesssim t^{-1/2}\|f\|_{L^1}.$
Since $L^1(\R)$ is dense in $H^1(\R)$, it follows that
$\|e^{-itL}(Qw(0)+w_1)\|_{L^\infty}\to0$ as $t\to\infty,$ and that
\begin{equation}
  \label{eq:nyan1}
\begin{split}
& \|Qw(t)\|_{L^\infty} \\ \le
&\|Qw(t)-e^{-itL}(Qw(0)+w_1)\|_{H^1}
+\|Qe^{-itL}(Qw(0)+w_1)\|_{L^\infty}
\\ \to & 0 \quad\text{as $t\to\infty$.}
\end{split}  
\end{equation}
Analogously to \eqref{eq:3.6}, we have
\begin{equation}
  \label{eq:nyan2}
\|Pw(t)\|_{L^\infty}\lesssim |E(t)-E_*|\|Qw(t)\|_{L^\infty}.  
\end{equation}
Combining \eqref{eq:nyan1} and \eqref{eq:nyan2}, we have
$\lim_{t\to\infty}\|Pw(t)\|_{H^1}=0.$
Thus by \eqref{eq:1.3} and \eqref{eq:trans}, 
$$
\lim_{t\to\infty}\left\|u(t)-e^{-i\theta(t)}\phi_{E(t)}
{}-e^{-itL}Q(w(0)+w_1)\right\|_{H^1}=0.$$
{}From \cite{AY} and \cite{Wed}, 
we see that there exists a $v_+\in H^1$ such that
$$
\lim_{t\to\infty}\|e^{-itL}Q(w(0)+w_1)-We^{it\pd_x^2}v_+\|_{H^1}=0,
$$
where $W=\lim_{t\to\infty}e^{-itL}e^{-it\pd_x^2}$.
Thus we complete the proof of Theorem \ref{thm1}.
\end{proof}

\section{Linear estimates}
\label{sec:4}
Let $R(\lambda)=(\lambda-L)^{-1}$ and let  $dE_{ac}(\lambda)$ be
the absolute continuous part of the spectral measure of $L$.
We have $R(\lambda-i0)=R(\lambda+i0)$ for $\lambda<0$ and
it follows from the spectral decomposition theorem that
\begin{equation*}
\begin{split}
Q e^{-itL}f=& \int_{-\infty}^\infty e^{-it\lambda}dE_{ac}(\lambda)f
\\ =& \frac{1}{2\pi i}\int_{-\infty}^\infty e^{-it\lambda}
Q(R(\lambda-i0)-R(\lambda+i0))fd\lambda.
\end{split}  
\end{equation*}
To prove Lemmas \ref{lem:2.2}--\ref{lem:2.4}, we will apply
Plancherel's theorem to the above formula. 
\subsection{High energy estimate}
To begin with, we will estimate the high frequency part of
the resolvent operators $R(\lambda\pm i0)$.
Let $\chi(x)$ be a smooth function satisfying
$0\le \chi(x)\le 1$  for $x\in\R$ and
$$\chi(x)=
\begin{cases}
1 \quad & \text{if $x\ge 2$,}\\
0 \quad & \text{if $x\le 1$,}
\end{cases}$$
and let $\chi_M(x)$ be an even function satisfying
$\chi_M(x)=\chi(x-M)$ for $x\ge0$.
\begin{lmm}
  \label{lem:4.4}
Assume (V1) and (V2). Then there exist positive numbers $M$ and $C$ such that 
\begin{align*}
& \sup_x\|\chi_M(\mathstrut{\sqrt{\lambda}})R(\lambda\pm i0)u\|_{L^2_\lambda(0,\infty)}
\le C\|u\|_{L^2(\R)},  \\
& \sup_x\|\chi_M(\mathstrut{\sqrt{\lambda}})\pd_x (R(\lambda-i0)-R(\lambda+i0))
u\|_{L^2_\lambda(0,\infty)}
\le C\|u\|_{H^{1/2}(\R)}
\end{align*}
for every $u\in \mathcal{S}(\R)$.
\end{lmm}
\begin{lmm}
  \label{lem:4.5}
Assume (V1) and (V2). Then there exist positive numbers $M$ and $C$ such that
$$\sum_{i=0,1}\|\chi_M(\mathstrut{\sqrt{\lambda}})\pd_x^iR(\lambda\pm i0)u
\|_{L^\infty_{x,\lambda}(\R^2)} \le C\|u\|_{L^1_x(\R)}$$
for every $\lambda\in\R$ and $u\in \mathcal{S}(\R)$.
\end{lmm}
\begin{proof}[Proof of Lemmas \ref{lem:4.4} and \ref{lem:4.5}]
Let $R_0(\lambda)=(\lambda+\pd_x^2)^{-1}$ and 
$$G_1(x,k)=\frac{e^{ik|x|}}{2ik}, \quad G_2(x,k)=\frac{e^{-k|x|}}{-2k}
.$$
We remark that $R_0(\lambda\mp i0)\delta=G_1(x,\pm k)$
for $\lambda=k^2$ with $k\ge0$ and $R_0(\lambda)\delta= G_2(x,k)$
for $\lambda=-k^2$ with $k>0$. 
If $M$ is sufficiently large, we have
\begin{equation}
  \label{eq:bo1}
R(\lambda\pm i0)u
=\sum_{j=0}^\infty R_0(\lambda\pm i0)(V R_0(\lambda\pm i0))^j u
\end{equation}
for $\lambda\in\R$ with $|\lambda|\ge M$ and $u\in \mathcal{S}(\R)$ since
$$\|\la \cdot\ra^{-1}R_0(\lambda\pm i0)\la \cdot\ra^{-1}\|_{B(L^2(\R))}
\lesssim \la \lambda\ra^{-1/2}.$$
By the definition of $G_1(x,k)$ and  Plancherel's theorem, 
\begin{equation}
  \label{eq:g1}
  \begin{split}
& \sup_x\int_\R dk \la k\ra 
 \left|\chi_M(k)(G_1(\cdot,k)*u)(x)\right|^2
\\  \lesssim &
\sup_x\int_\R dk\la k\ra^{-1} \left(
\left|\int_x^\infty u(y)e^{-iky}dy\right|^2
+\left|\int^x_{-\infty} u(y)e^{iky}dy\right|^2\right)dk
\\ \lesssim &
\|u\|_{L^2_x}^2.    
  \end{split}
\end{equation}

For $k=\mathstrut{\sqrt{\lambda}}\ge0$, it holds that
\begin{equation}
  \label{eq:comp}
  \begin{split}
& F_{1,n}(x,\pm k):= R_0(\lambda \mp i0)(V R_0(\lambda \mp i0))^nu(x)
\\ =&
\int_{\R^{n+1}}G_1(x-x_1,\pm k)\prod_{j=1}^n
\left(V(x_j)G_1(x_j-x_{j+1},\pm k)\right)
u(x_{n+1})dx_1\cdots dx_{n+1}.
  \end{split}
\end{equation}
Combining Minkowski's inequality with \eqref{eq:g1},
we have for $n\ge1$,
\begin{equation}
  \label{eq:bo2}
\begin{split}
  & \|\chi_M(\mathstrut{\sqrt{\lambda}})F_{1,n}(x,\pm\mathstrut{\sqrt{\lambda}})
\|_{L^2_\lambda(0,\infty)}
\\ \lesssim &
\int_{\R^{n+1}}dx_1\cdots dx_n \prod_{j=0}^{n-1}\{|V(x_{j+1})|
\sup_{|k|\ge M}(|kG_1(x_j-x_{j+1},k)|)\}
\\ & \quad \times
\left\{\int_{\R}dk \chi_M(k)^2\la k\ra^{-2n+1}
 \left|(G_1(\cdot,k)*u)(x_n)\right|^2 \right\}^{1/2}
\\ \lesssim & \|V\|_{L^1}^n \sup_{x_n}
\left(\int_\R dk k^{-2n+1}\chi_M(k)^2
\left|(G_1(\cdot,k)*u)(x_n)\right|^2 \right)^{1/2}
\\ \lesssim & M^{-n+1/2} \|V\|_{L^1}^n \|u\|_{L^2},
\end{split}  
\end{equation}
where $x_0=x$. Similarly, we have
\begin{equation}
  \label{eq:bo3}
 \|\chi_M(\mathstrut{\sqrt{\lambda}})\pd_xF_{1,n}(x,\pm\mathstrut{\sqrt{\lambda}})
\|_{L^2_\lambda(0,\infty)}
\lesssim  M^{-n+3/2} \|V\|_{L^1}^n \|u\|_{L^2}.
\end{equation}
Since 
\begin{align*}
\pd_x(R_0(\lambda-i0)-R_0(\lambda+i0))u=&\frac{1}{2ik}\pd_x\int_\R dy
u(y)(e^{-ik(x-y)}+e^{ik(x-y)})dy 
\\=& \sqrt{\frac{\pi}{2}}\left(e^{ikx}\hat{u}(k)-e^{-ikx}\hat{u}(-k)\right),
\end{align*}
it follows from  Plancherel's identity
\begin{equation}
  \label{eq:bo4}
  \begin{split}
\|\pd_x(R_0(\lambda-i0)-R_0(\lambda+i0))u\|_{L^2_\lambda}  \lesssim &
\left(\int_\R dk \la k\ra |\hat{u}(k)|^2dk\right)^{1/2}
\\ \lesssim & \|u\|_{H^{1/2}}.
   \end{split}
\end{equation}
Combining \eqref{eq:bo1}, \eqref{eq:bo2}--\eqref{eq:bo4}, we obtain Lemma \ref{lem:4.4}.
\par
Next, we will prove Lemma \ref{lem:4.5}.
In view of \eqref{eq:comp}, we have
\begin{equation}
  \label{eq:4.5.1}
  \begin{split}
& \sup_{x,k}\left(|\chi_M(k)F_{1,n}(x,k)|+|\chi_M(k)\pd_xF_{1,n}(x,k)|\right)
\\ \lesssim &
\sup_{x\in\R ,|k|\ge M} \la k\ra^{-n}\int_{\R^{n+1}}\prod_{j=1}^n|V(x_j)|
|u(x_{n+1})|dx_1\cdots dx_{n+1}
\\ \lesssim & M^{-n}\|V\|_{L^1(\R)}^n\|u\|_{L^1(\R)}.
  \end{split}
\end{equation}
For $\lambda=-k^2$ with $k>0$, we have
\begin{equation*}
  \begin{split}
& F_{2,n}(x,k):= R_0(\lambda)V R_0(\lambda)^nu(x)
\\ =&
\int_{\R^{n+1}}G_2(x-x_1,k)\prod_{j=1}^n \left(V(x_j)G_2(x_j-x_{j+1},k)\right)
u(x_{n+1})dx_1\cdots dx_{n+1},
  \end{split}
\end{equation*}
and it follows that
\begin{equation}
  \label{eq:4.5.2}
  \begin{split}
& \sup_{x,k}\left(|\chi_M(k)F_{2,n}(x,k)|+|\chi_M(k)\pd_xF_{2,n}(x,k)|
\right) \\ \lesssim & M^{-n}\|V\|_{L^1(\R)}^n\|u\|_{L^1(\R)}.
  \end{split}
\end{equation}
Combining \eqref{eq:4.5.1} and \eqref{eq:4.5.2}, we obtain
Lemma \ref{lem:4.5}.
\end{proof}
\par

\subsection{Low energy estimate}
Next, we will estimate the low frequency part of  $R(\lambda\pm i0)$.
Let $\widetilde{\chi}_M(x)=1-\chi_M(x)$.
\begin{lmm}
 \label{lem:5.1}
Assume (V1) and (V2). Let $M$ be a positive number given in
Lemma \ref{lem:4.4}.  Then there exists a positive number $C$ such that
for every $u\in \mathcal{S}(\R)$,
\begin{gather*}
\sup_x\|\la x\ra^{-3/2} \widetilde{\chi}_M(\mathstrut{\sqrt{\lambda}})
R(\lambda\pm i0)u\|_{L^2_\lambda(0,\infty)} \le C\|u\|_{L^2(\R)},
\\
\sup_x\|\widetilde{\chi}_M(\mathstrut{\sqrt{\lambda}}) \pd_x R(\lambda\pm i0)u
\|_{L^2_\lambda(0,\infty)} \le C\|u\|_{L^2(\R)}.
\end{gather*}
\end{lmm}

\begin{lmm}
  \label{lem:5.2}
Assume (V1) and (V2). Let $M$ be a positive number given in Lemma
\ref{lem:4.5}. Then there exists a positive number $C$ such that
$$
\sum_{j=0,1}\sup_{\lambda\in[-M,M]}
\left\| \la x\ra^{-1} \pd_x^jR(\lambda\pm i0)u\right\|_{L^\infty_x}
 \le C\|\la x\ra u\|_{L^1_x(\R)}$$
for every $\lambda\in\R$ and $u\in \mathcal{S}(\R)$.
\par
Furthermore, if $\sup_{x\in\R}e^{\alpha |x|}|V(x)|<\infty$ holds
for an $\alpha>0$, there exists a positive number $C$ such that
$$\sup_{\lambda\in[-M,M]} \left\|
\pd_x R(\lambda\pm i0)u\right\|_{L^\infty_x}\le C\|u\|_{L^1_x(\R)}
$$
for every $\lambda\in\R$ and $u\in \mathcal{S}(\R)$.
\end{lmm}

Before we start to prove Lemmas  \ref{lem:5.1} and  \ref{lem:5.2},
we recall some properties of the Jost functions. 
We refer the readers to Deift-Trubowitz \cite{DT} for the details. 
Let $f_1(x,k)$ and $f_2(x,k)$ be the solutions to $ Lu=k^2 u$ satisfying
$$
 \lim_{x\to\infty}|e^{-ikx}f_1(x,k)-1|=0,\quad
\lim_{x\to-\infty}|e^{ikx}f_2(x,k)-1|=0,
$$
and let $m_1(x,k)=e^{-ikx}f_1(x,k)$ and $m_2(x,k)=e^{ikx}f_2(x,k)$.
For each $x$, $m_1(x,k)$ and $m_2(x,k)$ are analytic in $k$ with $\Im k>0$,
continuous in $k$ with $\Im k\ge0$,  and satisfy
\begin{gather*}
  m_1(x,k)=1+\int_x^\infty \frac{e^{2ik(y-x)}-1}{2ik}V(y)m_1(y,k)dy,
\\ m_2(x,k)=1+\int^x_{-\infty} \frac{e^{2ik(x-y)}-1}{2ik}V(y)m_2(y,k)dy.
\end{gather*}
Deift-Trubowitz \cite{DT} tells us that for $x\in \R$ and $k\in\C$ with
$\Im k\ge0$,
\begin{gather}
  \label{eq:m1}
\left|m_1(x,k)-1\right|\lesssim 
\la k\ra^{-1}(1+\max(-x,0))\int_x^\infty dy \la y\ra|V(y)|,
\\
\label{eq:m2}
\left|m_2(x,k)-1\right|\lesssim 
\la k\ra^{-1}(1+\max(x,0))\int^x_{-\infty} dy \la y\ra|V(y)|,
\\
  \label{eq:m1'}
\left|\pd_xm_1(x,k)\right|\lesssim 
\la k\ra^{-1}\int_x^\infty dy \la y\ra|V(y)|,
\\
\label{eq:m2'}
\left|\pd_x m_2(x,k)\right|\lesssim 
\la k\ra^{-1}\int^x_{-\infty} dy \la y\ra|V(y)|.
\end{gather}
 For every $\delta>0$,
there exists a $C_\delta>0$ such that for $x\in \R$ and $k\in\C$ with $\Im k\ge0$ and $|k|\ge \delta$,
\begin{gather}
  \label{eq:m1"}
\left|m_1(x,k)-1\right|\le C_\delta\int_x^\infty dy |V(y)|,
\\
\label{eq:m2"}
\left|m_2(x,k)-1\right|\le C_\delta\int^x_{-\infty} dy|V(y)|.
\end{gather}

There exist continuous functions
$T(k)$, $R_1(k)$ and $R_2(k)$ on $\R$ satisfying
\begin{align}
\label{eq:alg1}
 & f_2(x,k)=\frac{R_1(k)}{T(k)}f_1(x,k)+\frac{1}{T(k)}f_1(x,-k),
\\ \label{eq:alg2}
& f_1(x,k)=\frac{R_2(k)}{T(k)}f_2(x,k)+\frac{1}{T(k)}f_2(x,-k)
\end{align}
for $k\in\R$ with $k\ne0$.
Let $[f(x),g(x)]=f'(x)g(x)-f(x)g'(x)$ and let $W(k)=[f_1(x,k),f_2(x,k)]$.
As is well known, the Wronskian $W(k)$ does not depend on $x$ and
$W(k)=2ik/T(k)\ne0$ for $k\ne0$. Moreover, the assumption (V2) implies
 $W(0)\ne0$.

\begin{proof}[Proof of Lemma \ref{lem:5.1}]
For $\lambda=k^2$ with $k\ge0$, the resolvent operator
 $R(\lambda\pm i0)$ has the kernel 
\begin{equation}
  \label{eq:kernel}
K_\pm(x,y,k)=\left\{
\begin{aligned}
   & -\frac{f_1(x,\pm k)f_2(y,\pm k)}{W(\pm k)} \quad \text{for $x>y$.,}
\\ & -\frac{f_2(x,\pm k)f_1(y,\pm k)}{W(\pm k)} \quad \text{for $x<y$.}
\end{aligned}\right.  
\end{equation}
That is,
\begin{equation}
  \label{eq:green}
\begin{split}
R(\lambda\pm i0)u=& -\frac{f_1(x,\pm k)}{W(\pm k)}\int_{-\infty}^x
dyf_2(y,\pm k)u(y)
\\ &  -\frac{f_2(x,\pm k)}{W(\pm k)}\int_x^\infty dyf_1(y,\pm k)u(y)
\\=:& I(\pm k)+II(\pm k).
\end{split}  
\end{equation}
\par

We will estimate $L^2_\lambda$-norm of the right-hand side of
\eqref{eq:green}. We may assume $x>0$.  Let
$$I=-\frac{f_1(x,k)}{W(k)}(I_1+I_2+I_3),$$
where
\begin{align*}
& I_1=\int_0^x dy f_2(y,k)u(y),\\
& I_2=\int_{-\infty}^0 dy e^{-iky}u(y)=\sqrt{2\pi}
\mathcal{F}_y(\chi_{(-\infty,0]}u)(k),\\
& I_3=\int_{-\infty}^0 dy e^{-iky}(m_2(y,k)-1)u(y).
\end{align*}
By \eqref{eq:m1} and \eqref{eq:m2}, we see
\begin{equation}
  \label{eq:5.1.1}
\sup_{x>0}\left(|f_1(x,k)|+\la x\ra^{-1}|f_2(x,k)|\right)<\infty,  
\end{equation}
\begin{equation}
  \label{eq:5.1.2}
\begin{split}
|I_1| \lesssim & \int_0^x dy \la y\ra|u(y)| \lesssim \la x\ra^{3/2}\|u\|_{L^2},
\end{split}  
\end{equation}
and
\begin{equation}
  \label{eq:5.1.3}
\sup_{x>0}|I_3|\lesssim 
\left\| \int_{-\infty}^x dy V(y)\right\|_{L^2_x(-\infty,0)}\|u\|_{L^2}
 \lesssim   \|u\|_{L^2}.
\end{equation}
Similarly, we have $II=-\frac{f_2(x,k)}{W(k)}(II_1+II_2)$ with
\begin{align*}
& II_1=\int_x^\infty dy e^{iky}u(y)=\sqrt{2\pi}
\mathcal{F}^{-1}(\chi_{[x,\infty)}u)(k),\\
& II_2=\int_x^\infty dy e^{iky}(m_1(y,k)-1)u(y),
\end{align*}
and
\begin{equation}
  \label{eq:5.1.4}
\sup_x|II_2|\lesssim \|u\|_{L^2}.  
\end{equation}
Obviously,
\begin{equation}
  \label{eq:5.1.5}
\sup_{x>0}(\|I_2\|_{L^2_k}+\|II_1\|_{L^2_k})\lesssim \|u\|_{L^2}.
\end{equation}
Since $W(k)\ne0$ for every $k\in\R$ and $\widetilde{\chi}_M(k)$ is compactly
supported, it follows from  \eqref{eq:5.1.1}--\eqref{eq:5.1.5} that 
$$
\sup_x\int_\R dk |k| \left|\widetilde{\chi}_M(k)
\int_\R dy K_\pm(x,y,k) u(y)\right|^2
\lesssim \la x\ra^3 \|u\|_{L^2}^2.$$
\par

By \eqref{eq:kernel}, we have
\begin{equation*}
\begin{split}
\pd_xR(\lambda\pm i0)u=& -\frac{\pd_xf_1(x,\pm k)}{W(\pm k)}
\int_{-\infty}^x dy f_2(y,\pm k)u(y)
\\ &  -\frac{\pd_x f_2(x,\pm k)}{W(\pm k)}
\int_x^\infty dyf_1(y,\pm k)u(y)
\\=& III(\pm k)+IV(\pm k).
\end{split}  
\end{equation*}
By symmetry, it suffices to consider the case where $x>0$. Let us rewrite
$III$ as 
$$III=-\frac{e^{ikx}}{W(k)}(III_1+III_2+III_3),$$
where
\begin{align*}
& III_1=(ik m_1(x,k)+\pd_x m_1(x,k))
\int_{-\infty}^0 dy f_2(y,k)u(y),\\
& III_2=ik m_1(x,k)\int_0^x dy f_2(y,k)u(y),\\
& III_3=\pd_x m_1(x,k)\int_0^x dy f_2(y,k)u(y).
\end{align*}
By \eqref{eq:m1} and \eqref{eq:m1'},
$$\sup_{k\in[-M,M]}\sup_{x\ge0}
\left(|km_1(x,k)|+|\pd_xm_1(x,k)|\right)<\infty.$$
Thus we have
\begin{equation}
  \label{eq:5.1.6a}
    \|III_1\|_{L^2_k(-M,M)+L^\infty_k(-M,M)}\lesssim \|u\|_{L^2}
\end{equation}
in the same way as \eqref{eq:5.1.3} and \eqref{eq:5.1.5}.
Under the assumption (V1) and (V2), we have
$1/T(k)\simeq k^{-1}$ as $k\to0$ and $R_1(k)$ and $R_2(k)$ are
continuous in $k\in\R$.
Hence by using \eqref{eq:alg1}, we see that
\begin{equation}
  \label{eq:5.1.6}
\|III_2\|_{L^2_k(-M,M)+L^\infty_k(-M,M)} \lesssim \|u\|_{L^2}
\end{equation}
follows in the same way as \eqref{eq:5.1.3} and \eqref{eq:5.1.5}.
By (V1), \eqref{eq:m1'} and Schwarz's inequality, we have
\begin{equation}
  \label{eq:5.1.7}
  \begin{split}
  |III_3|  \lesssim & \left(\int_x^\infty dy |V(y)|\right)
\int_0^x dy \la y\ra u(y)
\\ \lesssim & \la x\ra^{-1/2}\|\la x\ra^2 V\|_{L^1} \|u\|_{L^2}.
  \end{split}
\end{equation}
Similarly, we have
\begin{equation}
  \label{eq:5.1.8}
\sup_{x>0}\|IV\|_{L^2_k(-M,M)+L^\infty_k(-M,M)}\lesssim \|u\|_{L^2}.
\end{equation}
Combining \eqref{eq:5.1.6}--\eqref{eq:5.1.8},
we obtain
$$
\sup_{x>0}\int_\R dk |k| \left|\widetilde{\chi}_M(k)
\int_\R dy \pd_xK_\pm(x,y,k) u(y)\right|^2
\lesssim \|u\|_{L^2}^2.$$
Thus we complete the proof of Lemma \ref{lem:5.1}.
\end{proof}

\begin{proof}[Proof of Lemma \ref{lem:5.2}]
Since $W(k)$ is continuous and $W(k)\ne 0$ on $\R$, it follows from
\eqref{eq:m1}--\eqref{eq:m2'} and \eqref{eq:kernel} that
\begin{equation}
  \label{eq:ker-est}
\sup_{k\in[-M,M]}\sup_{x,y\in\R}
\la x\ra^{-1}|\pd_x^jK_\pm(x,y,k)|\la y\ra^{-1}<\infty
\quad\text{for $j=0$, $1$.}
\end{equation}
Thus we have 
\begin{equation}
  \label{eq:2.4.1}
\sup_{\lambda\ge0} \|\la x\ra^{-1} \widetilde{\chi}_M(\sqrt{\lambda})
\pd_x^jR(\lambda\pm i0)u\|_{L^\infty_x}\le C\|\la x\ra u\|_{L^1_x(\R)}
\quad\text{for $j=0$, $1$.}
\end{equation}
\par

For $\lambda<0$, the resolvent operator $R(\lambda)$ has the kernel
\begin{equation}
  \label{eq:kernel'}
K(x,y,\lambda)=\left\{
\begin{aligned}
  & -\frac{f_1(x,ik)f_2(y,ik)}{W(ik)}
 \quad \text{for $x>y$.,}\\
  & -\frac{f_2(x,ik)f_1(y,ik)}{W(ik)}
    \quad \text{for $x<y$,}
\end{aligned}\right.  
\end{equation}
where $k=\sqrt{-\lambda}$.
We have $W(ik)=0$ for a $k>0$ if and only if $\lambda=-k^2$
is an eigenvalue of $L$. Thus the assumption (V2) yields that
$W(ik)$ has a simple pole at $k=\sqrt{|E_*|}$ and
$W(ik)\ne0$ for $k\in[0,\sqrt{|E_*|})\cup(\sqrt{|E_*|},\infty)$.
Thus by \eqref{eq:m1}--\eqref{eq:m2'}, we have
$$\sup_{\lambda<0}\left(
\left\|\la x\ra^{-1}R(\lambda)Qu \right\|_{L^\infty_x}
+\left\|\la x\ra^{-1}\pd_xR(\lambda)Qu \right\|_{L^\infty_x}\right)
\lesssim \|\la x\ra u\|_{L^1_x}.
$$
Combining the above, we obtain the former part of Lemma \ref{lem:5.2}.
\par

Next, we will estimate $\pd_xK_\pm(x,y,k)$ and $\pd_xK(x,y,\lambda)$
assuming that $V(x)$ decays like $e^{-\alpha |x|}$. 
In view of \eqref{eq:m1}--\eqref{eq:m2'}, we have
\begin{equation}
  \label{eq:5.2.1}
\sup_{k\in[-M,M]}\left(\sup_{x>0>y}|\pd_xK_\pm(x,y,k)|+
\sup_{x<0<y}|\pd_xK_\pm(x,y,k)|\right)<\infty.  
\end{equation}
Suppose $x$ and $y$ has the same sign. By symmetry, we may assume $x>y>0$.
By \eqref{eq:m1}--\eqref{eq:m1'}, we have
$\sup_{k\in[-M,M]}|\pd_xm_2(y,k)|\lesssim \la y\ra \lesssim \la x\ra,$
and
\begin{align*}
|\pd_xm_1(x,k)m_2(y,k)|\lesssim &
\la x\ra\int_x^\infty dy \la y\ra |V(y)|
 \lesssim \la x\ra e^{-\alpha |x|}.
\end{align*}
As in the proof of Lemma \ref{lem:5.1}, it follows from  \eqref{eq:alg1} that
$$\sup_{k\in[-M,M]}\sup_{y>0}|km_2(y,k)|<\infty.$$
Combining the above, we see that
$$
\pd_xK_\pm(x,y,k)=\frac{-e^{\pm ik(x-y)}}{W(\pm k)}\{
\pm ikm_1(x,\pm k)+\pd_xm_1(x,\pm k)\}m_2(y,\pm k).$$
are uniformly bounded with respect to $x>y>0$ and $k\in[-M,M]$.
\par

By  \eqref{eq:m1}--\eqref{eq:m2"}, we have
\begin{gather*}
\sup_{\lambda\le -\alpha^2/16}\sup_{x,y\in\R}
|(\lambda-E_*)\pd_xK(x,y,\lambda)|<\infty,
\\ 
\sup_{\lambda<0}|\lambda-E_*|
\left(\sup_{x<0<y}|\pd_x K(x,y,\lambda)|
+\sup_{y<0<x}|\pd_x K(x,y,\lambda)|\right)<\infty.
\end{gather*}
Now we will prove the remaining case by using \eqref{eq:alg1} and
\eqref{eq:alg2}.
We may assume $x>y>0$ by symmetry.
Under the assumption $\sup_{x\in\R}e^{\alpha|x|}|V(x)|<\infty$,
$m_1(x,k)$ and $m_2(x,k)$ are
analytic in $k$ with $\Im k>-\alpha/2$ and there exists a $C_a>0$
for every $a>0$ such that
\begin{equation}
  \label{eq:m13}
|m_1(x,k)-1|\le C_a\int_x^\infty dy \la y\ra e^{-2\Im k y}|V(y)|
\end{equation}
for $x>-a$ and $-\alpha/2<\Im k<0$, and 
\begin{equation}
 \label{eq:m23}
 |m_2(x,k)-1|\le C_a\int^x_{-\infty} dy \la y\ra e^{2\Im k y}|V(y)|  
\end{equation}for $x<a$ and $-\alpha/2<\Im k<0$.
Furthermore, we see that
\begin{align*}
\frac{1}{T(k)}=&\frac{1}{2ik}[f_1(x,k),f_2(x,k)],\\
\frac{R_1(k)}{T(k)}=& \frac{1}{2ik}[f_2(x,k),f_1(x,-k)],\quad
\frac{R_2(k)}{T(k)}= \frac{1}{2ik}[f_2(x,-k),f_1(x,k)],\\
\end{align*}
are meromorphic in $k$ with $|\Im k|<\alpha/2$ and have a pole of order
$1$ at the origin.
Hence it follows from \eqref{eq:alg1}, \eqref{eq:m1} and \eqref{eq:m13} that
$$
\sup_{0<k\le \alpha/4}\sup_{y>0}k|m_2(y,ik)|<\infty.
$$
By \eqref{eq:m2} and \eqref{eq:m1'},
\begin{align*}
& \sup_{0<k\le \alpha/4}|m_2(y,ik)|\lesssim \la y\ra\lesssim \la x\ra, \\
& |\pd_xm_1(x,ik)|\lesssim \int_x^\infty dy \la y\ra|V(y)|
\lesssim \la x\ra e^{-\alpha x}.
\end{align*}
Hence, we have
\begin{align*}
& \sup_{0<k\le \alpha/4}\sup_{x>y>0} |\pd_xf_1(x,ik)f_2(y,ik)|
\\ =&
\sup_{0<k\le \alpha/4}\sup_{x>y>0} e^{-kx}
\left|(-km_1(x,ik)+\pd_xm_1(x,ik)) m_2(y,ik)\right| <  \infty.  
\end{align*}
Combining the above, we obtain
$\sup_{\lambda<0}\|\pd_xR(\lambda)Q\|_{B(L^1,L^\infty)}<\infty.$
Thus we prove the latter part of Lemma \ref{lem:5.2}.
\end{proof}

\subsection{Proof of Lemmas \ref{lem:2.2}--\ref{lem:2.4}}
Now, we are in position to prove Lemmas \ref{lem:2.2}--\ref{lem:2.4}.
\begin{proof}[Proof of Lemma \ref{lem:2.2}]
By the spectral decomposition theorem, we have
$$
Qe^{-itL}f=e^{-itL}\chi_M(L)f+Qe^{-itL}\widetilde{\chi}_M(L)f,
$$
and 
\begin{align}
\label{eq:lem4-1}
\chi_M(L)e^{-itL}f=& 
\frac{1}{2\pi i}\int_{-\infty}^\infty e^{-it\lambda}
\chi_M(\lambda)(R(\lambda-i0)-R(\lambda+i0))f d\lambda,
\\
Qe^{-itL}\widetilde{\chi}_M(L)f=&
\frac{1}{2\pi i}\int_{-\infty}^\infty e^{-it\lambda}
\widetilde{\chi}_M(\lambda)Q(R(\lambda-i0)-R(\lambda+i0))f d\lambda.
\end{align}  
Integrating \eqref{eq:lem4-1} by part, we see that
\begin{align*}
&\chi_M(L)e^{-itL}f \\ =&\frac{(it)^{-j}}{2\pi i}
\int_{-\infty}^\infty d\lambda e^{-it\lambda}
\pd_\lambda^j\{(R(\lambda+i0)-R(\lambda-i0))\chi_M(\lambda)\}f
 \quad\text{in $\mathcal{S}'_x(\R)$}   
\end{align*}
for any $t\ne0$ and $f\in \mathcal{S}_x(\R^2)$. Since
$$\|\pd_\lambda^jQR(\lambda\pm i0)
\|_{B(L^{2,(j+1)/2+0},L^{2,-(j+1)/2-0})}
\lesssim \la \lambda\ra^{-(j+1)/2},$$
the above integral absolutely converges in $L^{2,-(j+1)/2-0}_x$ for
$j\ge2$. 
\par

Suppose $g(t,x)=g_1(t)g_2(x)$,
$g_1\in C_0^\infty(\R\setminus\{0\})$ and $g_2\in\mathcal{S}(\R)$.
We define $\la \cdot,\cdot\ra_x$ and $\la \cdot,\cdot\ra_{t,x}$ as
$$\la u_1,u_2\ra_x:=\int_{-\infty}^\infty u_1(x)u_2(x)dx,\quad
\la v_1,v_2\ra_{t,x}:=\int_{-\infty}^\infty\int_{-\infty}^\infty
v_1(t,x)v_2(t,x)dxdt.$$

Making use of Fubini's theorem and integration by parts,
we have for $j\ge2$,
\begin{align*}
& \la \chi_M(L)e^{-itL}f,g\ra_{t,x}
\\=&
\frac1{2\pi i}\int_{-\infty}^\infty dt (it)^{-j}g_1(t)
\int_{-\infty}^\infty d\lambda e^{-it\lambda}
\pd_\lambda^j\left\la\chi_M(\lambda)(R(\lambda+i0)-R(\lambda-i0))f,
g_2\right\ra_x
\\=&
\frac{1}{2\pi i}\int_{-\infty}^\infty d\lambda
\pd_\lambda^j \left\la \chi_M(\lambda)(R(\lambda+i0)-R(\lambda-i0))f,
g_2 \right\ra_x
\int_{-\infty}^\infty dt (it)^{-j}g_1(t) e^{-it\lambda}
\\=& \frac{1}{\sqrt{2\pi}i}\int_{-\infty}^\infty d\lambda
(\mathcal{F}_tg_1)(\lambda)
\left\la \chi_M(\lambda)(R(\lambda+i0)-R(\lambda-i0))f,g_2\right\ra_x.
\end{align*}
Hence it follows from the above and Fubini's theorem that
\begin{align*}
& \left\la \chi_M(L)e^{-itL}f,g\right\ra_{t,x}
\\ = & \frac{1}{\sqrt{2\pi}i}
\int_{-\infty}^\infty dx \int_{-\infty}^\infty d\lambda 
\left(\chi_M(\lambda)(R(\lambda+i0)-R(\lambda-i0))f\right)
\mathcal{F}_t g(\lambda,x)
\end{align*}
for every $g\in C_0^\infty(\R_t\setminus\{0\})\otimes\mathcal{S}(\R_x)$.
Using Plancherel's theorem, we have
\begin{equation}
  \label{eq:pl1}
  \begin{split}
& \left|\la \chi_M(L)e^{-itL}f,g\ra_{t,x}\right|
\\ \le & (2\pi)^{-1/2}
\|\chi_M(\lambda)(R(\lambda+i0)-R(\lambda-i0))f\|_{L^\infty_xL^2_\lambda}
\|\mathcal{F}_t g(\lambda,\cdot) \|_{L^1_xL^2_\lambda}
\\ = & (2\pi)^{-1/2}
\|\chi_M(\lambda)(R(\lambda+i0)-R(\lambda-i0))f
\|_{L_x^\infty L^2_\lambda(0,\infty)}
\|g\|_{L_x^1L_t^2}.    
  \end{split}
\end{equation}
Similarly, we have
\begin{equation}
  \label{eq:pl1b}
\begin{split}
& |\la Qe^{-itL}\widetilde{\chi}_M(L)f \ra_{t,x}|
\\\le & (2\pi)^{-1/2} \|\la x\ra^{-3/2}\widetilde{\chi}_M(\lambda)
Q(R(\lambda-i0)-R(\lambda+i0))f\|_{L^\infty_xL^2_\lambda}
\|\la x\ra^{3/2}g\|_{L_x^1L_t^2},
\end{split}
\end{equation}
and
\begin{equation}
  \label{eq:pl1c}
\begin{split}
  & |\la \pd_x e^{-itL}Qf, g\ra_{t,x}| 
\\ \le & (2\pi)^{-1/2}\bigl(
\|\chi_M(\lambda)\pd_x(R(\lambda-i0)-R(\lambda+i0))f
\|_{L^\infty_xL^2_\lambda}
\\  & \phantom{(2\pi)^{-1/2}\bigl(}
+\|\widetilde{\chi}_M(\lambda)\pd_x(R(\lambda-i0)-R(\lambda+i0))Qf
\|_{L^\infty_xL^2_\lambda}\bigr)\|g\|_{L^\infty_xL^2_t}.    
\end{split}
\end{equation}
Since $C_0^\infty(\R_t\setminus\{0\})\otimes\mathcal{S}(\R_x)$
is dense in $L_x^1L_t^2$, Eqs. \eqref{eq:2.2.1} and \eqref{eq:2.2.2}
follow from \eqref{eq:pl1}--\eqref{eq:pl1c} and Lemmas \ref{lem:4.4}
and \ref{lem:5.1}.
By using the duality argument, we see that \eqref{eq:2.2.3} follows from
\eqref{eq:2.2.1}. Thus we complete the proof of Lemma \ref{lem:2.2}.
\end{proof}

To prove Lemma \ref{lem:2.3}, we need the following.
\begin{lmm}
  \label{lem:4.2} Assume (V1) and (V2).
Let $g(t,x)\in\mathcal{S}_\otimes(\R^2)$ and
$$U(t,x)=\frac{i}{\sqrt{2\pi}}\int_{-\infty}^\infty d\lambda
e^{-it\lambda}
\{R(\lambda-i0)+R(\lambda+i0)\}Q(\mathcal{F}^{-1}_tg)(\lambda,\cdot).
$$
Then,
\begin{align*}
U(t,x)=& 2\int_0^t ds e^{-i(t-s)L}Qg(s,\cdot)
+\int_{-\infty}^0 ds e^{-i(t-s)L}Qg(s,\cdot)
\\ & -\int_0^\infty ds e^{-i(t-s)L}Qg(s,\cdot).
\end{align*}
\end{lmm}
\begin{proof}
We may assume that $g(t,x)$ is written as $g(t,x)=g_1(t)g_2(x)$ with
$g_1$, $g_2\in \mathcal{S}(\R)$.
Let $h\in \mathcal{S}(\R)$ and 
\begin{align*}
& f(\lambda)=\la Q\{R(\lambda-i0)+R(\lambda+i0)\}g_2,h\ra,
\\ & f_\varepsilon(\lambda)=
\la Q\{R(\lambda-i\varepsilon)+R(\lambda+i\varepsilon)\}g_2,h\ra.
\end{align*}
Then $f(\lambda)$ and $f_\varepsilon(\lambda)$ are smooth
functions satisfying
$$\sup_{\lambda\in\R,\varepsilon>0}\la \lambda \ra^{k+1/2}
\left(\left|\pd_\lambda^kf(\lambda)\right|+
\left|\pd_\lambda^kf_\varepsilon(\lambda)\right|\right)<\infty$$
for every $k\in \N\cup \{0\}$ (see e.g \cite{Mu}), and 
\begin{equation}
  \label{eq:nyan3}
  \begin{split}
\int_\R U(t,x)\overline{h(x)}dx = & \frac{i}{\sqrt{2\pi}}
\int_\R d\lambda e^{-it\lambda}f(\lambda)(\mathcal{F}^{-1}g_1)(\lambda)
\\ =& \frac{i}{\sqrt{2\pi}} \int_\R ds \hat{f}(t-s)g_1(s).
\end{split}
\end{equation}
By the spectral decomposition theorem,
\begin{equation}
  \label{eq:fe} f_\varepsilon(\lambda)= \int_\R 
\frac{2(\lambda-\mu)}{(\lambda-\mu)^2+\varepsilon^2}d\la E_{ac}(\mu)g_2,h\ra.  
\end{equation}
Taking the Fourier transform of \eqref{eq:fe} and using Fubini's theorem,
we have
\begin{align*}
  \hat{f}_\varepsilon(t)= & \frac{1}{\sqrt{2\pi}} \int_\R d\la E_{ac}(\mu)g_2,h\ra
\int_\R d\lambda e^{-it\lambda}
\frac{2(\lambda-\mu)}{(\lambda-\mu)^2+\varepsilon^2}
\\=& \sqrt{2\pi}i \int_\R d\la E_{ac}(\mu)g_2,h\ra
e^{-it\mu-\varepsilon|t|}\operatorname{sgn}t.
\end{align*}
Hence it follows that
\begin{equation}
  \label{eq:hf}
\begin{split}
\hat{f}(t)= & \lim_{\varepsilon\downarrow0} \hat{f_\varepsilon}(t)
 = \sqrt{2\pi}i \sgn t \la Qe^{-itL}g_2,h\ra.
\end{split}  
\end{equation}
Substituting \eqref{eq:hf} into \eqref{eq:nyan3}, we obtain
\begin{align*}
\la U(t,\cdot), h\ra=& 
{}-\int_\R ds g_1(s)\left(\int_\R \sgn(t-s)e^{-i(t-s)\mu}d\la E_{ac}(\mu)g_2,h\ra
\right)
\\ =&
\int^t_{-\infty} ds g_1(s)\la e^{-i(t-s)L}g_2,h\ra
{}-\int_t^\infty ds g_1(s)\la e^{-i(t-s)L}g_2,h\ra.
\end{align*}
Thus we complete the proof of Lemma \ref{lem:4.2}.
\end{proof}

\begin{proof}[Proof of Lemma \ref{lem:2.3}]
Since $\mathcal{S}_\otimes(\R^2)$ is dense in $L_x^1L_t^2$, it suffices
to prove \eqref{eq:2.3.1} for $g\in \mathcal{S}_\otimes(\R^2)$.
As in the proof of Lemma \ref{lem:2.2}, we have
$$
\|\la x\ra^{-1}\pd_x^jU(\cdot,x)\|_{L_x^\infty L_t^2}
\le \|\la x\ra^{-1}\pd_x^j\{R(\lambda-i0)+R(\lambda+i0)\}Q
\mathcal{F}^{-1}_tg(\lambda,\cdot)\|_{L^\infty_xL^2_\lambda}.
$$
Applying Plancherel's theorem and Minkowski's inequality, we have
\begin{align*}
& \|\la x\ra^{-1}\pd_x^jU(\cdot,x)\|_{L_x^\infty L_t^2}
\\ \le & \left\|
\|\la \cdot\ra^{-1}\pd_x^j\{R(\lambda-i0)+R(\lambda+i0)\}Q
\la\cdot\ra^{-1}\|_{B(L^1_x,L^\infty_x)}
\|\la\cdot\ra\mathcal{F}^{-1}_tg(\lambda,\cdot)
\|_{L^1_x}\right\|_{L^2_\lambda}
\\ \le & \sup_\lambda 
\|\la \cdot\ra^{-1}\pd_x^j \{R(\lambda-i0)+R(\lambda+i0)\}Q
\la \cdot\ra^{-1}\|_{B(L^1_x,L^\infty_x)}\|\la x\ra g\|_{L_x^1 L_t^2}
\end{align*}
for $j=0$, $1$.
Hence it follows from Lemmas \ref{lem:4.5} and \ref{lem:5.2} that
\begin{equation}
  \label{eq:U-est}
\|\la x\ra^{-1}U\|_{L_x^\infty L_t^2}+\|\la x\ra^{-1}\pd_xU
\|_{L^\infty_xL^2_t} \lesssim \|\la x\ra g\|_{L_x^1 L_t^2}.
\end{equation}

For $I=[0,\infty)$ and $I=(-\infty,0]$, we have
\begin{align*}
& \int_I e^{-i(t-s)L} Qg(s)ds \\=&  \int_\R ds \chi_I(s)
\left(\int_\R e^{-i(t-s)\lambda} dE_{ac}(\lambda)g(s,\cdot)\right)
\\=&\frac{1}{2\pi i}\int_{\R^2}d\lambda ds
e^{-i(t-s)\lambda}\{R(\lambda-i0)-R(\lambda+i0)\}Q\chi_I(s)g(s,\cdot)
\\=& -i\mathcal{F}_\lambda\{(R(\lambda-i0)-R(\lambda+i0))Q
\mathcal{F}^{-1}_s(\chi_I(s)g)(\lambda,\cdot)\}(t).
\end{align*}
By Plancherel's identity and Minkowski's inequality, we have
\begin{equation}
  \label{eq:2.2.5}
  \begin{split} & 
\left\|\la x\ra^{-1}\pd_x^j \int_I e^{-i(t-s)L}Qg(s)ds
\right\|_{L^\infty_xL_t^2}
\\ \le & \left\|\la x\ra^{-1}\pd_x^j(R(\lambda-i0)-R(\lambda+i0))Q
\mathcal{F}^{-1}_s(\chi_I(s)g)(\lambda,\cdot)\right
\|_{L_x^\infty L_\lambda^2}
\\ \le &   \sup_\lambda
\|\la\cdot\ra^{-1}\pd_x^j \{R(\lambda-i0)-R(\lambda+i0)\}Q
\la \cdot\ra^{-1}\|_{B(L^1_x,L^\infty_x)}
\|\la x\ra g\|_{L_x^1 L_t^2}
  \end{split}
\end{equation}
for $j=0$, $1$. Combining \eqref{eq:U-est}--\eqref{eq:2.2.5} with Lemma
\ref{lem:4.2}, we obtain \eqref{eq:2.3.1}.
Since \eqref{eq:2.3.2} can be obtained in exactly the same way,
we omit the proof.
\end{proof}
Finally, we will prove Lemma \ref{lem:2.4}.
To prove Lemma \ref{lem:2.4}, we will use a lemma of Christ and Kiselev
\cite{ChK}.

\begin{proof}[Proof of Lemma \ref{lem:2.4}]
Let $(q,p)=(4,\infty)$ or $(q,p)=(\infty,2)$ and let
$$Tg(t)=\int_\R ds e^{-i(t-s)L}Qg(s).$$
Lemmas \ref{lem:2.1} and \ref{lem:2.2} imply
$f:=\int_\R ds e^{isL}Qg(s)\in L^2(\R)$ and
$$  \|Tg(t)\|_{L^q_tL^p_x}\lesssim  \|f\|_{L^2_x}
\lesssim  \|\la x\ra^{3/2}g\|_{L^1_xL^2_t}.$$
Thus by Schwarz's inequality,
we see that there exists a $C>0$ such that for every
$g\in \mathcal{S}(\R^2)$,
\begin{equation}
  \label{eq:T}
\|Tg(t)\|_{L^q_tL^p_x}\le C \|g\|_{L^2_tL^2_x(\R,\la x\ra^5dx)}.  
\end{equation}
Since $q>2$, it follows from Lemma 3.1 in \cite{SmS} and \eqref{eq:T}
that
\begin{equation}
  \label{eq:T1}
\left\|\int_{s<t} ds e^{-i(t-s)L}Qg(s)\right\|_{L_t^qL_x^p} \lesssim
\|g\|_{L^2_tL^2_x(\R,\la x\ra^5dx)}.  
\end{equation}
Thus we prove Lemma \ref{lem:2.4}.
\end{proof}

\section*{Acknowledgments} 
The author would like to express his gratitude to Professor Kenji
Nakanishi for telling the author the paper by Christ and Kiselev
\cite{ChK}.
This research is supported by Grant-in-Aid for Scientific
Research (No. 17740079).

\end{document}